\documentclass[11pt]{article}

\usepackage{latexsym}
\usepackage{color}
\usepackage{amsmath}
\definecolor{dred}{rgb}{.8,0,0}

\begin{document}
\newcommand{\qed}{\hphantom{.}\hfill $\Box$\medbreak}
\newtheorem{exa}{Example}[section]
\newtheorem{lemma}{Lemma}[section]
\newtheorem{remark}{Remark}[section]
\newtheorem{prop}{Proposition}[section]
\newtheorem{theo}{Theorem}[section]

\def \cA {{\cal A}}
\def \cB {{\cal B}}
\def \cC {{\cal C}}
\def \cD {{\cal D}}
\def \cE {{\cal E}}
\def \cG {{\cal G}}
\def \cI {{\cal I}}
\def \cL {{\cal L}}
\def \cK {{\cal K}}
\def  \cR  {{\cal R}} 
\def  \cS  {{\cal S}} 
\def \cT {{\cal T}}

\title{\bf Minimum embedding of any Steiner triple system into a $3$-sun system via matchings \thanks{G. Lo Faro and A. Tripodi were supported by INDAM (GNSAGA) and A. Tripodi was supported by FFABR Unime 2019}}

\author{Giovanni Lo Faro\\ \small
Dipartimento di Scienze Matematiche e Informatiche, \\  \small
Scienze Fisiche e Scien\-ze della Terra \\  \small
Universit\`a di Messina, Messina, Italia\\ \small email:
lofaro@unime.it, atripodi@unime.it
\and
 Antoinette Tripodi\\ \small
Dipartimento di Scienze Matematiche e Informatiche, \\  \small
Scienze Fisiche e Scien\-ze della Terra \\  \small
Universit\`a di Messina, Messina, Italia\\ \small email:
lofaro@unime.it, atripodi@unime.it}
\date{}
\maketitle

\begin{abstract}
Let $G$ be a simple finite graph and $G'$ be a subgraph of $G$. A $G'$-design $(X,\cB)$ of order $n$ is said to be {\em embedded} into a $G$-design $(X\cup U,\cC)$ of order $n+u$, if there is an injective function $f:\cB\rightarrow\cC$ such that $B$ is a subgraph of $f(B)$ for every $B\in\cB$. The function $f$ is called an \textit{embedding} of $(X,\cB)$ into $(X\cup U,\cC)$. If $u$ attains the minimum possible value, then $f$ is a \textit{minimum} embedding. Here, by means of K\"{o}nig's Line Coloring Theorem and edge coloring properties a complete solution is given to the problem of determining a minimum embedding of 
any $K_3$-design  (well-known
as Steiner Triple System or, shortly, STS) into a 3-sun system or, shortly, a 3SS (i.e., a $G$-design where $G$ is a graph on six vertices consisting of  a triangle with three pendant edges which form a 1-factor).
\end{abstract}



\section{Introduction} 

If $G$ is a graph, then let $V(G)$ and $\cE(G)$ denote the vertex-set and edge-set of $G$, respectively. 
Given a set $\Gamma$ of  pairwise non-ismorphic simple graphs, a {\em $\Gamma$-design of order $n$}   
 is  a pair $(X, {\cal B})$ where 
${\cal B}$ is a collection of graphs (called {\em blocks}) each  isomorphic to some element of $\Gamma$, whose edges partition 
$\cE(K_n)$, where $K_n$ is the complete graph  of order $n$ 
on $X$; 
  if the edges of the blocks of ${\cal B}$ partition a proper spanning subgraph of  $K_n$, then we speak of  {\em partial  $\Gamma$-design of order $n$}. 
  If $\Gamma=\{G\}$, then we simply write $G$-design. Let $\Sigma(G)$ denote the set of all integers $n$ such that there exists a $G$-design of order $n$. A $K_3$-design of order $n$ is known as Steiner triple  system and denoted by STS$(n)$; it is well-known that $\Sigma(K_3)=\{n\in N : n\equiv 1,3 \pmod{6} \}$.

Let $G$ be a simple finite graph and $G'$ be a subgraph of $G$. A $G'$-design $(X,\cB)$ of order $n$ is said to be {\em embedded} into a $G$-design $(X\cup U,\cC)$ of order $n+u$, if there is an injective function $f:\cB\rightarrow\cC$ such that $B$ is a subgraph of $f(B)$ for every $B\in\cB$. The function $f$ is called an \textit{embedding} of $(X,\cB)$ into $(X\cup U,\cC)$.  If $u$ attains the minimum possible value, then $f$ is a \textit{minimum} embedding. Note that a special case occurs when $G=G'$ and the related embedding problem is better known as Doyen-Wilson problem (see \cite{dw,FLLHH,lt1,lt2, LT}).

The embedding problems have interesting applications to networks (\cite{CQS}), that is why  
they have been investigated in several papers. 
In particular, the minimum embedding problem of STSs into $G$-designs have been studied in the case when $G=K_4$,  $G=K_4 - e$ (the complete graph on four vertices with one deleted edge), or $ G=K_3+ e$ (a \emph{kite}, i.e., a triangle with one pendant edge)  have been solved in 
\cite{CLQ}, \cite{CLQ2}, \cite{GQR}, \cite{MR}. 

In \cite{FJLLS} the authors embed a cyclic STS of order   $n\equiv 1 \pmod{6} $   into a  3-sun system of order $2n-1$,  i.e., a $G$-design  where $G$ is a graph on six vertices consisting of  a triangle with three pendant edges which form a 1-factor, and as an open problem they ask whether it is possible to embed any STS into a 3-sun system. Here we give an answer to this open problem by determining  the minimum embedding for any Steiner triple system.
More precisely, for every integer $n\in \Sigma (K_3)$ detoted by $u_{min} (n)$ the minimum integer $u$ such that any STS$(n)$ can be embedded into a 3-sun system of order $n+u$, as main result we  prove the following theorem.

\bigskip

\noindent \textbf{Main Theorem}
\begin{itemize} 
\item [(i)] If $n\equiv 1,3,9, 19 \pmod{24} $, then  $u_{min} (n)= \frac{n-1}{2}$ for every $n\neq 3,9$,  $u_{min} (3)=6$, and $u_{min} (9)=7$. 
\item [(ii)] If $n\equiv 7,13,15,21 \pmod{24} $, then $u_{min} (n)= \frac{n-1}{2}+2$ for every $n\neq 7,13$,  $u_{min} (7)=6$, and $u_{min} (13)=11$. 
\end{itemize}

To obtain our result we  make use of some results on edge colorings and, in particular, of  K\"{o}nig's Line Coloring Theorem, which here, for convenience,  is  formulated in terms of matchings (for 
definitions and results on edge colorings or matchings, the reader is referred to  \cite{BM}).

\begin{theo} \label{konig} \ {\rm (K\"{o}nig's Line Coloring Theorem)}
Let $G$ be a bipartite multigraph with maximum
degree $\Delta$. Then $\cE(G)$ can be partitioned into  $M_1, M_2, $ $\ldots,  M_{\Delta}$ such that each $M_i$,  $1 \leq i \leq  \Delta$, is a matching in $G$.
\end{theo}

\section{Notation and basic results}

In what follows, we will denote:
\begin{itemize}
\item  the triangle on the vertices $a$, $b$ and $c$ by $(a,b,c)$;
\item  the  kite consisting of the triangle  $(a,b,c)$ and the pendant edge $\{c,d\}$ by $(a,b,c; d)$;
\item  the  \emph{bull} graph consisting of the triangle  $(a,b,c)$ and the pendant edges $\{b,d\}$ and $\{c,e\}$ by $(a,b,c; d, e)$;
\item  the  3-sun consisting of the triangle  $(a,b,c)$ and the pendant edges $\{a,d\}$, $\{b,e\}$  and $\{c,f\}$ by $(a,b,c; d, e,f)$.
\end{itemize}
If $G$ is a kite, a bull, or a 3-sun, then its triangle  will be denoted by $t(G)$. 

\bigskip

In this section we will give the necessary condition for embedding a Steiner triple system  into a 3-sun system  and prove some useful results to get our main result.  From now on, if  $f$ is an embedding of $(X,\cT)$ into $(X\cup U,\cS)$, then $f(\cT)$ will be denoted by $\cS_\cT$.
Finally, we recall that a 3-sun system of order $n$, or shortly a 3SS$(n)$, exists if and only if $n\equiv 0,1,4,9 \pmod{12} $   (see \cite{FJLLS}).

\begin{lemma} \label{CN}
If there exists a $3$SS$(n+u)$ embedding an STS$(n)$, then  $u\geq \frac{n-1}{2}$.
\end{lemma}
{\em Proof.}  Since an  STS$(n)$ has $\frac{n(n-1)}{6}$ triples, then in order to complete every triple so to abtain a 3-sun, necessarily $n\cdot u\geq 3\frac{n(n-1)}{6}$ and so $u\geq \frac{n-1}{2}$.
\hfill$\Box$

\bigskip

 In general, to construct a $3$SS$(m)$ $(X\cup U,\cS)$
embedding  a  STS$(n)$  $(X,\cT)$,   we need to complete each triangle of $\cT$ to a 3-sun by using some edges of the complete bipartite graph $K_{n, m-n}$ on $X\cup U$ and partition into 3-suns the remaining edges of $K_{n, m-n}$ along with those of the complete graph $K_{n, m-n}$ on $U$. In the following lemma a partial $3$SS embedding an STS$(n)$ is constructed by using all the edges  of the above complete bipartite graph.

\bigskip

\begin{lemma} \label{partial}
Any STS$(n)$, $n\geq 7$, can be embedded into a partial $3$SS$(\frac{3n-1}{2})$.
\end{lemma}
{\em Proof.}  Let $(X,\cT)$  be an STS$(n)$ and consider its incidence graph $\cI$, i.e., the bipartite graph 
 whose vertex set is $X\cup \cT$ and whose edges are determined by joining $x\in X$ to $t\in \cT$ if and only if $x\in t$. 
 In the graph $\cI$ every vertex of $X$ has degree $\frac{n-1}{2}$ and every vertex of $ \cT$ has degree 3. Since the maximum degree of $\cI$ is $\Delta=\frac{n-1}{2}$, by K\"{o}nig's Line Coloring Theorem  the edges of $\cI$  can be partitioned into $\Delta$ matchings $M_1, M_2, \ldots, M_{\Delta}$, each of which satures the vertices of $X$, i.e., every vertex of  $X$ is incident to an edge of each matching. 
Let $\cal S$ be the set of $3$-suns on  $X\cup \{M_1, M_2, \ldots, M_{\Delta}\}$ obtained by completing each triple of $\cT$  to a 3-sun as follows: for every $t=(x_1, x_2,x_3)\in \cT$, consider the 3-sun  $(x_1, x_2,x_3; M_{i_1}, M_{i_2}, M_{i_3})$, where $\{x_j ,t \}\in M_{i_j} $ for every $j=1,2,3$. $(X\cup \{M_1, M_2, \ldots, M_{\Delta}\}, \cal S)$ is  a partial  $3$SS$(\frac{3n-1}{2})$ embedding $(X,\cT)$.
\hfill$\Box$

\bigskip

The lower bound given by Lemma \ref{CN} is attained if $n\equiv 1,3,9, 19 \pmod{24} $, $n \neq 3,9$, as it is established by the following proposition.

\begin{prop} \label{u_min}
For every $n\equiv 1,3,9, 19 \pmod{24} $, $n\geq 19$, $u_{min} (n)= \frac{n-1}{2}$.
\end{prop}
{\em Proof.}  Let $(X,\cT)$ be  any STS$(n)$ with $n\equiv 1,3,9, 19 \pmod{24} $, $n\geq 19$. By Lemma \ref{partial}, it can be embedded into a partial $3$SS$(\frac{3n-1}{2})$ $(X\cup \{M_1, M_2, \ldots, M_\frac{n-1}{2}\}, \cS)$. Since $\frac{n-1}{2}\equiv 0,1,4, 9 \pmod{12} \geq 9$,  there exists a 3SS$(\frac{n-1}{2})$  $(\{M_1, M_2, \ldots, M_\frac{n-1}{2}\}, \cS')$. Then $(X\cup\{M_1, M_2, \ldots, M_\frac{n-1}{2}\}, \cS\cup \cS' )$ is a 3SS$(\frac{3n-1}{2})$ which embeds $(X,\cT)$. 
\hfill$\Box$


\begin{lemma}
If $n=3,9$, then  $u_{min} (n)=6, 7$, respectively.
\end{lemma}
{\em Proof.}  Any STS$(3)$ can be trivially embedded into a 3SS of any admissible order $v \geq 9$ and so $u_{min} (3)=6$. 

Let $(X\cup U,\cS)$ be a 3SS$(9+u)$ {embedding} an STS$(9)$ $(X,\cT)$. By Lemma \ref{CN} $u\geq 4$. If $u=4$, then $\cS\setminus \cS_\cT$ contains only one 3-sun $S$ such that $V(S)\subseteq U$, which is impossible and so $u_{min} (9)\geq 7$. 
To prove that $u_{min} (9)= 7$,   on $Z_{16}$ we give  the blocks of a 3SS embedding the unique STS$(9)$ (whose triangles are in bold):
\bigskip

$({\mathbf{0,1,2}}; 9,10,11),\ 
(({\mathbf{0,3,6}};10,15,9),\ 
({\mathbf{0,4,8}};11,9,13),$

$({\mathbf{0,5,7}};12,9,15),\
({\mathbf{1,3,8}};9,10,11),\ 
({\mathbf{1,4,7}};11,10,9),$

$({\mathbf{1,5,6}};12,10,13),\ 
({\mathbf{2,3,7}};9,11,12),\
({\mathbf{2,4,6}};10,11,12),$

$({\mathbf{2,5,8}};12,13,14),\ 
({\mathbf{3,4,5}};9,12,14),\ 
({\mathbf{6,7,8}};11,10,9),$

$(0,13,15;14,7,8),\ 
(1,14,15;13,4,9),\ 
(3,12,14;13,15,11)$

$(2,13,14;15,4,7),\
(5,11,12;15,7,10),\ 
(6,10,14;15,8,9),$

$(9,12,13;10,8,11),\ 
(10,11,15;13,9,4).$
 \hfill$\Box$

\begin{lemma} \label{CN2}
Let $n\equiv 7,13,15,21 \pmod{24} $. If there exists  a $3SS(n+u)$ embedding an STS$(n)$, then  $u\geq \frac{n-1}{2}+2$. 
\end{lemma}
{\em Proof.}  Let $n=24k+r$, $r \in \{7,13,15,21\}$. If  $(X, \cS)$ is a 3SS$(n+u)$ embedding an STS$(n)$, then by Lemma \ref{CN} $n+u\geq \frac{3n-1}{2}=36k+ \frac{3r-1}{2}$, where $\frac{3r-1}{2}\in \{10,19,22,31\}$. Since $n+u\equiv 0,1,4,9 \pmod{12} $,  this implies $u\geq \frac{n-1}{2}+2$ 
\hfill$\Box$

\bigskip

\begin{remark} \label{remark1}
{\rm For every $n\equiv 7,13,15,21\!\!\! \pmod{24} $, if $(X\cup U, \cS)$  is a $3$SS$(n+\frac{n+3}{2})$ embedding an STS$(n)$ $(X, \cT)$,   
then  each vertex $x\in X$ appears 
in exactly two block of $\cS \setminus \cS_\cT$ as a pendant vertex (therefore, for every $S\in \cS\setminus \cS_\cT$ the vertices of $t(S)$ are in $U$).}
\end{remark}

\bigskip

The lower bound established by Lemma \ref{CN2} is not attained when $n=7,13$, as it is showed by the following lemma.

\begin{lemma}
If $n=7,13$, then  $u_{min} (n)=6, 11$, respectively.
\end{lemma}
{\em Proof.}  Let $(X\cup U,\cS)$ be a 3SS$(n+u)$ {embedding} an STS$(n)$ $(X,\cT)$, where $n=7,13$. By Lemma \ref{CN2}, $u \geq \frac{n+3}{2}.$ \\
If $n=7$  and $u=5$, then $|\cS\setminus \cS_\cT|= 4$, whereas by Remark \ref{remark1} $|\cS\setminus \cS_\cT|> 4$, and so $u_{min} (7)\geq 6$. To prove that $u_{min} (7)= 6$,   on $Z_{13}$ we give  the blocks of a 3SS embedding the unique STS$(7)$:
\bigskip

$({\mathbf{0,1,2}}; 7,8,9),\ 
(({\mathbf{0,3,4}};8,7,9),\ 
({\mathbf{0,5,6}};9,8,10),\
({\mathbf{1,3,5}}; 9,8,7),$

$({\mathbf{1,4,6}}; 10,7,12),\ 
({\mathbf{2,3,6}}; 7,9,8),\ 
({\mathbf{2,4,5}};8,11,9),$

$(0,10,11;12,5,3),\ 
(1,7,12;11,8,5), \ 
(2,10,12;11,9,8),\
(4,8,10;12,9,3),$

$(9,11,12;7,5,3),\ 
(6,7,11;9,10,8).$
\bigskip

\noindent If $n=13$  and $u=8$, then $|\cS\setminus \cS_\cT|= 9$. Therefore,  by Remark \ref{remark1} a partial triple system on $U$ with 9 triangles should be exist, which is impossible because a maximun packing of $K_8$ with triangles (i.e., a partial $K_3$-design of order $8$ with the maximum number of blocks) have 8 blocks, and so  $u_{min} (13)\geq 11$. Since there are two non-isomorphic STS$(13)$s, in order to prove that $u_{min} (13)= 11$ we need to embed each STS$(13)$. Firstly, we embed the cyclic one  into a 3SS  on $Z_{24}$ as follows:
\bigskip

$({\mathbf{0,1,4}};13,18,14),\ 
({\mathbf{1,2,5}};13,23,14),\ 
({\mathbf{2,3,6}};13,18,14),$

$({\mathbf{3,4,7}};13,15,14),\ 
({\mathbf{4,5,8}};13,15,14),\ 
({\mathbf{5,6,9}};13,18,19),$

$({\mathbf{6,7,10}};13,15,14),\ 
({\mathbf{7,8,11}};13,15,14),\ 
({\mathbf{8,9,12}};13,20,15),$

$({\mathbf{9,10,0}};13,15,14),\ 
({\mathbf{10,11,1}};13,16,14),\ 
({\mathbf{11,12,2}};13,16,14),$

$({\mathbf{12,0,3}};13,15,23),\ 
({\mathbf{0,2,7}};16,22,21),\ 
({\mathbf{1,3,8}};15,19,16),$

$({\mathbf{2,4,9}};15,16,14),\ 
({\mathbf{3,5,10}};14,16,17),\ 
({\mathbf{4,6,11}};17,15,18),$

$({\mathbf{5,7,12}};17,16,14),\ 
({\mathbf{6,8,0}};16,17,18),\ 
({\mathbf{7,9,1}};17,15,16),$

$({\mathbf{8,10,2}};18,16,17),\ 
({\mathbf{9,11,3}};16,15,17),\ 
({\mathbf{10,12,4}};18,17,19),$

$({\mathbf{11,0,5}};17,19,18),\ 
({\mathbf{12,1,6}};18,17,19),$

$(0,17,20;21,6,1),\ 
(1,19,21;22,2,3),\ 
(2,16,18;20,3,4),$

$(3,15,20;22,13,4),\ 
(4,21,22;23,2,0),\ 
(5,19,20;21,7,6),$

$(6,21,23;22,8,0),\ 
(7,18,20;22,9,8),\ 
(8,19,22;23,10,5),$

$(9,17,21;22,13,10),\ 
(10,20,22;23,11,12),\ 
(11,19,23;21,12,1),$

$(12,20,21;23,13,14),\ 
(13,14,16;18,15,17),\ 
(13,22,23;19,11,5),$

$(14,17,19;18,15,16),\ 
(14,20,23;22,16,7),\ 
(15,16,21;19,23,13),$

$(15,18,22;23,19,16),\ 
(17,18,23;22,21,9).$

\bigskip

\noindent A 3SS$({24})$ embedding the non cyclic STS$(13)$ can be obtained from the above one by replacing the 3-suns 

\bigskip

$({\mathbf{0,1,4}};13,18,14),\ 
({\mathbf{0,2,7}};16,22,21),$
 
$({\mathbf{2,4,9}};15,16,14),\ 
({\mathbf{7,9,1}};17,15,16),$

\bigskip

\noindent with

$({\mathbf{9,1,4}};14,18,16),\ 
({\mathbf{9,2,7}};15,22,17),$
 
$({\mathbf{0,2,4}};16,15,14),\ 
({\mathbf{0,1, 7}};13,16,21),$
\hfill$\Box$
\bigskip

In order to prove that for every $n\equiv 7,13,15,21 \pmod{24} $, $n\neq 7,13$,  $u_{min} (n)$ equals the lower bound of Lemma \ref{CN2}, it will be useful the  following lemma.

\begin{lemma} \label{M+1}
\ {\rm (\cite{BM})} Let $M$ and $N$ be disjoint matchings of a graph $G$ with $|M| > |N|$. Then there are disjoint matchings $M'$ and $N'$ of $G$ such that $|M'| =|M|-1$,
$|N'| = |N|+1$ and $M' \cup N'= M\cup N$.
\end{lemma}



\bigskip

Now, we determine $u_{min} (n)$  for every $n\equiv 7,13, 15,  21 \pmod{24} $ with the exception of few small orders, which will be settled in Section \ref{small}.  

 In graph theory, the degree of a vertex of a 
graph is the number of edges that are incident to the vertex; here, we define \emph{$2$-degree} of a vertex $x$ of a 
$\Gamma$-design $\cD$,  and denote by $d_2(x)$, the number of blocks of ${\cal D}$ containing $x$ as a vertex of degree $2$. The \emph{$2$-degree sequence} of $\cD$ is the non-decreasing sequence of its vertex 2-degrees.  

In what follows,  if $G$ is a graph whose vertices belong to  $Z_u$, then we call \emph{orbit of $B$ under $Z_u$}  the set $(G)= \{G+i : i\in Z_u\}$, where $G+i$ is the graph with $V(G+i)=\{a+i : a\in V(G)\}$ and $\cE(G+i)=\{\{a+i , b+i\}: \{a, b\}\in V(G)\}$.

\begin{lemma} \label{bull}
For any $u=12k+h$, $h=5, 8, 9, 
12$  and $k\geq 3$, there exists a $\{$bull, $3$-sun$\}$-design of order ${u}$  whose $2$-degree sequence  is $(2,3,3,3,3,4,4,$ $\ldots, 4)$.
\end{lemma}
{\em Proof.}  Consider the following orbits under $Z_u$: for $i=1,2,3$,   $\cB_i=(B_i)$,   where $B_1=(0, 6k-2,  4k+3; 3k, 6k-1)$, $B_2=( 6k,0,  4k+1; 6k+2, 6k+1)$,  and $B_3=(0, 6k-1,  4k+2; 3k, 6k)$;  for 
$j=0,1,\ldots,k-4$, $\cS_j=(S_j)$,  where $S_j=(5k+1+j, 5k-j, 0; 3k, k, u-2-2j)$. On $Z_u$ define the set of graphs $\cA=(\cB_1 \cup \cB_2^\ast \cup \cB_3^\ast\cup \cB) \cup [ (\cup _{j=0} ^{k-4}\cS_j)  \cup  \cS^\ast\cup \cS] $, where $\cB_2^\ast =\cB_2\setminus \{B_2\}$, $\cB_3^\ast =\cB_3\setminus \{B_3+i: i=0,4k+1, 6k, 6k+1, 6k+2\}$, $\cS^\ast =\{(6k-1,  4k+2, 0; 3k, 6k, 4k+1),(10k,  8k+3, 4k+1; 7k+1, 10k+1, 6k),(12k-1,  10k+2, 6k; 9k, 12k, 0),(12k,  10k+3, 6k+1; 9k+1, 12k+1, 4k+1),(12k+1,  10k+4, 6k+2; 9k+2, 12k+2, 0)\}$, while $\cB $ and $ \cS$ depend on $h$.\\
\emph{$a)$ $h=5$}: \ 
$\cB $ is the orbit of $(6k+1, 0,  3k; 3k+2, 6k+3)$ under $Z_u$; $\cS=\emptyset$.\\
\emph{$b)$ $h=8$}: \ 
$\cB = \{(6k+3+i,i,  3k+i; 6k+4+i, 9k+1+i), (9k+5+i,3k+2+i,  6k+2+i; 6k+4+i, 12k+3+i): i=0,1,\ldots,3k+1,\  i\in Z_u\}\cup \{(12k+7+i, 6k+4+i, 9k+4+i; 9k+5+i, 3k-3+i): i=0,1,\ldots,6k+3,\  i\in Z_u\}$;  
$\cS=\{(i, 3k+2+i, 9k+6+i; 3k+1+i, 6k+3+i, 6k+4+i)): i=0,1,\ldots,3k+1,\  i\in Z_u\}$.\\
\emph{$c)$ $h=9$}: \ 
$\cB$ is the orbit of $(6k+1, 0,  3k; 3k+3, 9k+3)$ under $Z_u$;  
$\cS=\{(3i, 3k+2+3i, 6k+4+3i;6k+5+3i, 9k+7+3i, 9k+6+3i)): i=0,1,\ldots,4k+2,\  i\in Z_u\}$.\\
\emph{$d)$ $h=12$}: \ 
$\cB = \{(6k+1+i, i, 3k+i; 6k+6+i, 9k+5+i), (9k+4+i, 3k+3+i,  6k+3+i; 6k+6+i, 12k+8+i): i=0,1,\ldots,3k+2,\  i\in Z_u\}\cup \{(12k+7+i,6k+6+i,  9k+6+i; 12k+9+i, 3k-1+i): i=0,1,\ldots,6k+5,\  i\in Z_u\}$;  
$\cS=\{(i, 3k+3+i, 9k+9+i; 6k+3+i, 9k+6+i, 6k+6+i)): i=0,1,\ldots,3k+2,\  i\in Z_u\}\cup \{(3i, 3k+2+3i, 6k+4+3i;6k+8+3i, 9k+10+3i, 9k+6+3i)): i=0,1,\ldots,4k+3,\  i\in Z_u\}$.\\
 $(Z_u, \cA)$ is the required design, where $d_2(6k)=2$, the vertices $d_2(0)= d_2(4k+1)=d_2(6k+1)=d_2(6k+2)=3$, and the  remaining vertices have 2-degree 4.
 \hfill$\Box$

\begin{prop} \label{u_min2}
For every $n\equiv 7,13, 15,
 21 \pmod{24} $, $n\geq 79$, $u_{min} (n)= \frac{n+3}{2}$.
\end{prop}
{\em Proof.}  Let $(X,\cT)$  be an STS$(n)$, $n\equiv 7,13, 15, 
21 \pmod{24} $, $n\geq 79$,  and $\cI$ be its incidence graph. $\cE (\cI)$  can be partitioned into $\Delta=\frac{n-1}{2}$ matchings $M_1, M_2, \ldots, M_{\Delta}$ (see proof of Lemma \ref{partial}). 
By applying Lemma \ref{M+1} and by using similar arguments as the proof of Theorem 6.3 in \cite{BM}, it is possible to partition $\cE (\cI)$  into $\Delta+{2}$ matchings $M'_1, M'_2, \ldots, M'_{\Delta +2}$, such that  $M'_i$ 
covers the vertices of $X\setminus X_i$, where $|X_1|=2$,  $|X_i|=3$ for $i=2,3,4,5$, and $|X_i|=4$ for $ i =6,7,\ldots, {\Delta +2} $ (note that each vertex of $X$ is missing in exctaly two matchings). 
If $\cal S$ denotes the set of $3$-suns on  $X\cup \{M'_1, M'_2, \ldots, M'_{\Delta+2}\}$ obtained by completing each triple of $\cT$  as in the proof of Lemma \ref{partial}, the pair  $(X\cup \{M'_1, M'_2, \ldots, M'_{\Delta+2}\}, \cal S)$ is a partial  $3$SS$(\frac{3(n+1)}{2})$  embedding $(X,\cT)$. In order to complete the proof it will be sufficient to decompose the graph $K_{\Delta+2} \cup \cal M$ into 3-suns, where $K_{\Delta+2} $ is the complete graph based on $\{M'_1, M'_2, \ldots, M'_{\Delta+2}\}$ and  $\cal M$ is the bipartite graph on $X\cup \{M'_1, M'_2, \ldots, M'_{\Delta+2}\}$ such that  $\{x, M'_i\}\in \cE(\cal M)$ if and only if $x\in X$ is missing in $M'_i$. By using Lemma \ref{bull}, the complete graph $K_{\Delta+2} $ can be decomposed into bulls or 3-suns so that $d_2(M'_1)=2$, $d_2(M'_i)=3$ for  $i=2,3,4,5$, and $d_2(M'_i)=4$ for $ i =6,7,\ldots, {\Delta +2} $.  To obtain the required decompostion it is sufficient to complete each bull to a 3-sun using the edges of $\cal M$. \hfill$\Box$

\section{Cases left}\label{small}

To determine $u_{min} (n)$ for the  remaining orders    $n\in \{15,21,31,37,39,45,55,$ $61,63,69\}$, we will start from an STS$(n)$  $(X,\cT)$,  with $X=\{x_1,x_2, \ldots,x_n\}$, and prove that $(X,\cT)$ can be embedded in a 3-sun system 
 $(X \cup Z_{ \frac{n+3}{2}},\cS)$ by taking the following steps.\\
\emph{Step} $1.\ $ Partition the edges of 
 the complete graph on $Z_{ \frac{n+3}{2}}$ into a set ${\cal A}$ of triangles, kites, bulls or 3-suns so that $|{\cal A}|=|{\cal S}\setminus {\cal T}|=({n^2+20n+3)}/{48}$ and $\sum_{i=0} ^{(n+1)/2} d_2 (i)=2n$. For later convenience (see \emph{Step} $4.$), give ${\cal A}$  partitioned into suitable subsets ${\cal A}_j$, $j\in J$, such that for every  $j\in J$ and for every vertex $i\in Z_{ \frac{n+3}{2}}$, the number of blocks of ${\cal A}_j$ containing $i$ as a vertex of degree 2 is at most 1.\\
\emph{Step} $2.\  $ Partition the edge-set of the incidence graph $\cI$ of $(X,\cT)$ into $\frac{n+3}{2}$  matchings $M_0, M_1, \ldots, M_{\frac{n+1}{2}}$ such that, denoted by $X_i$   the set of vertices of $X$ not satured by $M_i$,  $|X_i|= d_2 (i)$ for each $ i =0,1,\ldots, {\frac{n+1}{2}}$. \\
\emph{Step} $3.\ $ Complete each triple of $\cT$  as in the  proof of Lemma \ref{partial} and obtain a partial 3-sun system  
$(X\cup \{M_0, M_1, \ldots, M_{\frac{n+1}{2}}\}, \cal S)$   
embedding $(X,\cT)$.\\  
\emph{Step} $4.\ $ Call  \emph{missing graph} 
 the bipartite graph  $\cal M$ on $X\cup \{M_0, M_1,$ $ \ldots, $ $M_{\frac{n+1}{2}}\}$ consisting of all the edges $\{x, M_i\}$ such that  $x\in X_i$ and, for the sake of simplicity, for every $ i =0,1,\ldots, {\frac{n+1}{2}}$ identify $M_i$ with $i\in Z_{ \frac{n+3}{2}}$. \\
\emph{Step} $5.\ $  Partition the edges of  the missing graph into suitable matchings $M'_j$, $j\in J$, such that for  every $j\in J$ the edges of  $M'_j$ can be used to complete the blocks of  ${\cal A}_j$  so to obtain a 3-sun system of order $\frac{3(n+1)}{2}$ embedding  $(X,\cT)$.

\bigskip

To begin with, we give an alternative solution for  $n\equiv 15 \pmod{24} $ (which settles the orders $v=15, 39, 63$ as well) by means of a technique used in \cite{FJLLS} and involving the concepts of parallel classes and resolution of an STS. 

A  \emph{parallel class} of an STS$(n)$ is a set of  $\frac{n}{3}$ triples such that no two triples in the set share an element;  a partition of all triples of an STS$(n)$  into parallel classes is a \emph{resolution} and the STS is said to be \emph{resolvable}. An STS$(n)$ together with a resolution of its triples is a \emph{Kirkman triple system}, KTS$(n)$, and exists if and only if $n\equiv 3 \pmod{6} $(see \cite{HAND}).

\begin{prop} \label{15}
For every $n\equiv 15 \pmod{24} $, $u_{min} (n)= \frac{n+3}{2}$.
\end{prop}
{\em Proof.}  Let $(X,\cT)$  be an STS$(n)$,  $n=24k+15$, $k\geq 0$. 
 Consider a resolution $P_i$, $ i=1,2,\ldots, 6k+4$  of a KTS  on $Z_{ \frac{n+3}{2}}$. Without loss of generality, assume that $P_1$ contains the triangle $t=(0,1,2)$. Construct a  set  $\cal K$ of  kites obtained by attaching the edges of $t$ to the  triangles $t_1, t_2, t_3$ of $P_2$  containing  $0, 1, 2 $, respectively, and  the set ${\cal A}_0$ of 3-suns   obtained from the parallel classes $P_i$, $ i=5,6,\ldots, 6k+4$ by using the technique in Lemma 3.8 of \cite{FJLLS}. The set ${\cal A}=\cup_{j=0}^{4}{\cal A}_j $, where  ${\cal A}_1=P_1 \setminus \{t\}$, ${\cal A}_2=(P_2 \setminus \{t_1,t_2,t_3\})\cup \cK $ and ${\cal A}_j =P_j$ for $j=3,4$, is a partition of 
$\cE(K_{\frac{n+3}{2}})$ such that $\sum_{i=0} ^{(n+1)/2} d_2 (i)=2n$. After applying \emph{Step} $2.\  $, \emph{Step} $3.\ $  and \emph{Step} $4.\ $ proceed as follows. 
It is easy to see that the missing graph admits two matchings $M'_1$ and $M'_2$ both saturing the vertices  $3, 4,  \ldots,  {\frac{n+1}{2}}$; while, the edges of $\cal M$  not in $M'_1$ and $M'_2$ form a subgraph  with maximun degree 2 and so can be partitioned into two  matchings $M'_3$ and $M'_{4}$ both saturing all the vertices  of $Z_{ \frac{n+3}{2}}$. For every $j=1,2,3,4$, complete the blocks of  ${\cal A}_j$ by using the edges of $M'_j$.
\hfill$\Box$



\begin{prop} \label{smallorders}
For every $n\in \{21,31,37,45,55, 61,69\}$, $u_{min} (n)= \frac{n+3}{2}$.
\end{prop}
{\em Proof.} Let $(X,\cT)$  be an STS$(n)$. 

For $n=21$, partition the edges of the complete graph on $Z_{12}$ into the following set ${\cal A}$: 
\begin{tabbing}
\ \ \ \  ${\cal A}_1=$ \= $\{(1,2,0;11), (3,7,2;5), (0,4,3;9)\}\}$ \\
\ \ \ \  ${\cal A}_2=$\>$\{(0,5,6),(1,8,11), (7,4,10), (2,9,8;10), (3,1,5;10,8)\}$ \\
\ \ \ \  ${\cal A}_3=$ \> $\{(0,9,10), (3,6,8), (5,7,11), (2,4,11;6), (1,7,9;6,11)\}$ \\
\ \ \ \  ${\cal A}_4=$ \>$\{ (0,7,8), (3,10,11), (5,9,4;8), (1,4,6;9), (2,6,10;5)\}$ 
\end{tabbing}
\noindent where $d_2(i)=3$
for $i\in \{ 5,6,8,9,10,11\}$ and $d_2(i)=4$  for $i\in \{ 0,1,2,3,4,7\}$. After applying \emph{Step} $2.\  $, \emph{Step} $3.\ $  and \emph{Step} $4.\ $ proceed as follows. Since $\cal M$ has maximun degree $4$, it is easy to see that $\cal M$ admits a matching $M'_1$ saturing $\{ 0,1,2,3,4,7\}$.  Use  $M'_1$ to complete the kites in ${\cal A}_1$. The graph obtained from  ${\cal M}$ by deleting the edges of $M'_1$ is a bipartite graph such that all the vertices in $Z_{12}$  has  degree 3 and so its edges can be partitioned into three matchings $M'_2$, $M'_3$ and $M'_{4}$, each of which satures the vertices  of $Z_{ 12}$. For every $j=2,3,4$, use the edges of $M'_j$ to complete the blocks of  ${\cal A}_j$. 

For $n=31$, partition the edges of the complete graph on $Z_{17}$ into the following set ${\cal A}$: 

\begin{tabbing}
\ \ \ \ ${\cal A}_1=$ \= $\{(0,4,1;7)+i: i=2,3,4,5,11,12,13,14,\  i\in Z_{17}\} \cup$\\
 ~\> $\{(10,12,0;3,7)\}$ \\
\ \ \ \  ${\cal A}_2=$\>$\{(0,4,1;7)+i: i=0,1,6,7,8,9,15,16,\  i\in Z_{17}\}\cup \{(14,7,9;2,0)\}$ \\
\ \ \ \  ${\cal A}_3=$ \> $\{(0,7,2;10)+i: i=1,4,13,15,16,\  i\in Z_{17}\}\cup\{(10,14,11;0), $\\
 ~\> $(9,4,2;12,0),(12,2,14;10,5)\}$ \\
\ \ \ \  ${\cal A}_4=$ \>$\{(0,7,2;10)+i: i=3,5,6,8,9,11,14,\  i\in Z_{17}\}$ 
\end{tabbing}

\noindent where $d_2(i)=2$
for $i\in \{ 0,2,7\}$ and $d_2(i)=4$  for $i\in Z_{17}\setminus  \{ 0,2,7\}$. After applying \emph{Step} $2.\  $, \emph{Step} $3.\ $  and \emph{Step} $4.\ $ proceed as follows. Consider a subgraph  $\cal M'$ of the missing graph such that each vertex in $Z_{17}$ has degree 2. Partition  the edges of $\cal M'$ into two matchings $M'_1$ and $M'_2$ and use them  to  complete the kites in ${\cal A}_1$ and ${\cal A}_2$, respectively. After deleting the edges of $M'_1$ and $M'_2$  the remaining edges of  ${\cal M}$ can be partitioned into two matchings $M'_3$ and $M'_{4}$, each of which satures the vertices  in $Z_{17}\setminus  \{ 0,2,7\}$ and can be used to  complete the kites in ${\cal A}_3$ and ${\cal A}_4$, respectively. 

By similar arguments it is possible to settle  the remaining cases $n=37,45,55, 61,69$, for which we refer to Appendix where we give the sets ${\cal A}_js$, which automatically determine the matchings $M'_js$. 
\hfill$\Box$

\section{Main result and conclusion}

Combining Lemmas \ref{CN}, \ref{CN2},  and Propositions \ref{u_min}, \ref{u_min2}, \ref{15}, \ref{smallorders}
gives our main result.

\bigskip

\noindent \textbf{Main Theorem}
\begin{itemize} 
\item [(i)] If $n\equiv 1,3,9, 19 \pmod{24} $, then  $u_{min} (n)= \frac{n-1}{2}$ for every $n\neq 3,9$,  $u_{min} (3)=6$, and $u_{min} (9)=7$. 
\item [(ii)] If $n\equiv 7,13,15,21 \pmod{24} $, then $u_{min} (n)= \frac{n-1}{2}+2$ for every $n\neq 7,13$,  $u_{min} (7)=6$, and $u_{min} (13)=11$. 
\end{itemize}

In \cite{LT}   a complete solution to the Doyen-Wilson problem for  3-sun systems is given  and it is proved that any $3$SS$(n)$ can be embedded in a $3$SS$(m)$ if and only if $m\geq \frac 75 n+1$ or $m=n$. For every  integer $v\in \Sigma (K_3)$,  combining Main Theorem with the above result gives an integer $m_v$ such that  any  STS$(v)$ can be embedded in a $3$SS$(m)$ for every admissible $m\geq m_v$.
A question to be asked is the following.

\bigskip

\noindent \textbf{Open Problem} 
 Can one embed any STS$(v)$ in a $3$SS$(m)$ for every admissible $m$ such that $ v + u_{min} (v) < m < m_v$?

\bigskip\bigskip\bigskip

\noindent {\Large \textbf{Appendix}}

\bigskip

\noindent {$n=37$}:
\begin{tabbing} 
${\cal A}_1=$ \= $\{(4,11,0;8)+2i: i=0,1,\ldots, 9,\  i\in Z_{20}\} $\\
${\cal A}_2=$ \> $\{(5,12,1;9)+2i : i=0,1,\ldots, 9,\  i\in Z_{20}\} $\\ ${\cal A}_3=$ \> $\{(14,16,13;2,19),(4,6,3;16,9)\}$\\ 
${\cal A}_4=$ \> $\{(1,3,0;6)+i : i=0,12,17,\  i\in Z_{20}\}\cup$\\  
~\> $\{(7,12,2;17), (16,17, 19;14),(5,7,4;17,10), (6,8,5;18,11),$\\  
~\> $(8,10,7;15,13), (11,9,8;4,14),(10,12,9;17,15),(2,4,1;14,16)\}$ \\ 
${\cal A}_5=$ \> $\{(1,3,0;6)+i : i=14,15,\  i\in Z_{20}\}\cup$\\  
~\> $\{(5,10,0;15), (8,13, 3;18),(0,2,19;4), (3,5,2;15,8),(7,9,6;19,12),$\\  
~\> $ (11,13,10;18,16),(12,14,11;19,17),(1,18,19;4,5),(6,11,1;16,7)\}$
\end{tabbing}

\bigskip

\noindent {$n=45$}:
\begin{tabbing} 
${\cal A}_1=$ \= $\{(1,13,7;19)+i: i=0,1,2,3,4,\  i\in Z_{24}\} \cup$\\
~\>  $\{(8,16,0;22,12),(9,17,1;23,19),(10,18,2;6,20),(19,11,3;0,21),$\\
~\>  $( 20,12,4;6,22),(21,13,5;19,23 ),( 22,14,6;20,0),( 7,15,23;21,12)\} $\\
${\cal A}_2=$ \= $\{(0,1,5)+3i: i=0,1,\ldots,7,\  i\in Z_{24}\} $\\
${\cal A}_3=$ \= $\{(1,2,6)+3i: i=0,1,\ldots,7,\  i\in Z_{24}\} $\\
${\cal A}_4=$ \= $\{(2,3,7)+3i: i=0,1,\ldots,7,\  i\in Z_{24}\} $\\
${\cal A}_5=$ \= $\{(1,3,10; 12,6,20)+i: i\in Z_{24}\setminus \{11,23\}\} \cup \{(0,2,9; 18,5,19), $\\
~\>  $(12,14,21;18,17,7)\}$
 \end{tabbing}
 
\bigskip

\noindent {$n=55$}:
\begin{tabbing} 
${\cal A}_1=$ \= $\{(13,27,0;25)+i: i=3,4,\ldots, 13,\  i\in Z_{29}\}\cup$\\  
~\> $\{(15,0,2;25,27), (12,14, 27;10,13),(28,13,15;0,11), (14,0,16;27,12),$\\  
${\cal A}_2=$ \> $\{(13,27,0;25)+i: i=1,17,18,\ldots, 28,\  i\in Z_{29}\} $\\ 
${\cal A}_3=$ \> $\{(0,10,11;2,6)+i:  i\in Z_{29}\}$\\  
${\cal A}_4=$ \> $\{(0,9,12;2,6)+i: i\in Z_{29}\}$
\end{tabbing}

\bigskip

\noindent {$n=61$}:
\begin{tabbing} 
${\cal A}_1=$ \= $\{(0,10,29;9)+i: i=0,1,2,3,6,17,18,19,20,\  i\in Z_{32}\}\cup$\\  
~\> $\{(14,22,6;30)\}\cup\{(23,10,13;22,29)+i: i=0,1,2,\  i\in Z_{32}\}\cup$\\  
~\> $\{(4,23,26;3,18)+i: i=0,1,3,4,5,\  i\in Z_{32}\}\cup\{(26,13,16;25,24), $\\  
~\> $(21,18,31;30,7)\}$\\  
${\cal A}_2=$ \> $\{(8,16,0;24)+i: i=0,1,\ldots,5,\  i\in Z_{32}\}\cup\{(4,14,1;13),$\\  
~\> $ (0,22,19;31),(2,24,21;1)\}\cup$\\  
~\> $\{(1,20,23;0,15)+i: i=0,2,5,\  i\in Z_{32}\}\cup\{(5,2,15;14,7)\}$\\  
${\cal A}_3=$ \> $\{(17,4,7;16,23)+i: i=0,1,2,3,4,5,\  i\in Z_{32}\}$\\  
${\cal A}_4=$ \> $\{(9,0,2;11,17)+i: i\in Z_{32}\}$\\  
${\cal A}_5=$ \> $\{(5,0,1;6,15)+i: i\in Z_{32}\}$
\end{tabbing}

\bigskip

\noindent {$n=69$}:

\begin{tabbing} 
${\cal A}_1=$ \= $\{(4,2,0;6,34)+3i: i=5,6,7,8,9,10,\  i\in Z_{36}\}$\\  
${\cal A}_2=$ \= $\{(4,2,0;6,34)+3i: i=0,1,2,3,4,11,\  i\in Z_{36}\}\cup$\\  
~\> $\{(24,12,0;30,18)+i, (30,18,6;24)+i: i=0,1,2,3,4,5,\  i\in Z_{36}\}$\\  
${\cal A}_3=$ \> $\{(0,7,15;1)+2i: i=0,1,\ldots,17,\  i\in Z_{36}\}$\\  
${\cal A}_4=$ \> $\{(1,8,16;2)+2i: i=0,1,\ldots,17,\  i\in Z_{36}\}$\\  
${\cal A}_5=$ \> $\{(9,20,0;3,13)+i:  i\in Z_{36}\}$\\  
${\cal A}_6=$ \> $\{(0,6,1;10,32,11)+9i,(1,7,2;5,11,12)+9i,(2,8,3;5,9,13)+9i,$\\ 
~\> $(3,9,4;6,14,8)+9i,(4,10,5;1,7,8)+9i,(5,11,6;35,8,9)+9i,$
\\ 
~\> $(6,12,7;16,9,17)+9i,(7,13,8;4,23,18)+9i : i=0,1,2,3,\  i\in Z_{36}\}$
\end{tabbing}


\bigskip




\begin{thebibliography}{66}

\newcommand{\bbt}{\bibitem}


\bbt{BM}  J.A. Bondy, U.S.R. Murty, Graph theory with applications, North Holland, 1976.


 \bibitem{HAND} C.J. Colbourn,  {\em Triple Systems}, in:  C. J. Colbourn and J. H. Dinitz (eds.), CRC Handbook of Combinatorial Designs, Second Edition, Chapman and Hall/CRC, Boca Raton, FL, 2007,  pp. 58--71.


\bbt{CLQ}  C.J. Colbourn, A.C.H. Ling, G. Quattrocchi, {\em Minimum embedding of Steiner triple systems into $(K_4-e)$-designs I},   Discrete Math., 308 (2008) 5308–-5311.

\bbt{CLQ2} C.J. Colbourn, A.C.H. Ling, G. Quattrocchi, {\em Minimum embedding of Steiner triple systems into $(K_4-e)$-designs II},   Discrete Math., 309 (2009) 400–-411.

\bbt{CQS} C.J. Colbourn, G. Quattrocchi, V.R. Syrotiuk, {\em Grooming for two-period optical networks},    Networks 52 (2008),  307–-324.

\bbt{dw} 
J. Doyen and R.M. Wilson, 
{\em Embeddings of Steiner triple systems}, 
Discrete Math. 5 (1973), 229-239.


\bbt{FJLLS} 
C.M. Fu,  N.H. Jhuang, Y.L. Lin, S.W. Lo, and H.M. Sung,
{\em From Steiner triple systems to $3$-sun  systems}, 
Taiwan. J. Math. 16  (2012), 531-543.

\bbt{FLLHH} 
C.M. Fu, Y.L. Lin, S.W. Lo, Y.F. Hsu, and W.C. Huang,
{\em The Doyen–Wilson theorem for bull designs}, 
Discrete Math. 313 (2013), 498-507.

\bbt{GQR} M. Gionfriddo, G. Quattrocchi, G. Ragusa, {\em Minimum embedding of STSs into $(K_3+e)$-systems},   Discrete Math., 313 (2013) 1419–-1428.

\bbt{K} D.  K\"{o}nig, {\em \"{U}ber Graphen und ihre Anwendung auf Determinantentheorie und Mengenlehre},  Math. Ann. 77(4) (1916),  453--465.


\bbt{lt1} 
G. Lo Faro and A. Tripodi, 
{\em Embeddings of $\lambda$-fold kite systems, $\lambda \geq 2$}, 
Australas. J. Combin. 36 (2006), 143-150.

\bbt{lt2} 
G. Lo Faro and A. Tripodi, 
{\em The Doyen-Wilson theorem for kite systems}, 
Discrete Math. 306 (2006), 2695-2701.


\bbt{LT} 
G. Lo Faro and A. Tripodi, 
{\em The Doyen-Wilson theorem for $3$-sun systems}, 
Ars Math. Contemp. 16 (2019), 119--139.

\bbt{MR} M. Meszka, A. Rosa, {\em Embedding Steiner triple systems into Steiner systems $S(2, 4, v)$},  Discrete Math., 274 (2004) 199–-212.





\end{thebibliography}
 \end{document}